\documentclass[12pt]{article}
\usepackage{amsfonts}
\usepackage{latexsym}
\usepackage{color, amsmath,amssymb, amsfonts, amstext,amsthm, latexsym, dsfont}
\usepackage{epic,eepic}
\usepackage[version=4]{mhchem}

\usepackage{longtable}
\usepackage[title]{appendix}
\usepackage{indentfirst}
\usepackage{graphicx}
\usepackage{epstopdf}
\usepackage{subfig}
\usepackage{float}

\usepackage[utf8]{inputenc}
\usepackage{url}
\usepackage[colorlinks=false, linktocpage=true]{hyperref}

\numberwithin{equation}{section}
\numberwithin{Lem}{section}
\numberwithin{Defi}{section}
\numberwithin{Theo}{section}
\numberwithin{Rem}{section}
\numberwithin{Coro}{section}
\numberwithin{Fig}{section}

\begin{document}

\title{Most Probable Transitions \\from Metastable to Oscillatory Regimes in a Carbon Cycle System}

\author{\bf\normalsize{
Wei Wei$^{1,}$\footnotemark[2],
Jianyu Hu$^{1,}$\footnotemark[1],
Jianyu Chen$^{1,}$\footnotemark[3]
and Jinqiao Duan$^{2,}$\footnotemark[4]
}\\[10pt]
\footnotesize{$^1$Center for Mathematical Sciences, Huazhong University of Science and Technology,} \\
\footnotesize{Wuhan, Hubei 430074, China.} \\[5pt]
\footnotesize{$^2$Departments of Applied Mathematics \& Physics, Illinois Institute of Technology, Chicago, IL 60616, USA.}
}

\footnotetext[2]{Email: \texttt{weiw16@hust.edu.cn}}
\footnotetext[1]{Email: \texttt{jianyuhu@hust.edu.cn}}
\footnotetext[3]{Email: \texttt{jianyuchen@hust.edu.cn}}
\footnotetext[4]{Email: \texttt{duan@iit.edu}}
\footnotetext[1]{is the corresponding author}

\date{}
\maketitle
\vspace{-0.3in}

\begin{abstract}
Global climate changes are related to the ocean's store of carbon. We study a carbonate system of the upper ocean, which has metastable and oscillatory regimes, under small random fluctuations. We calculate the most probable transition path via a geometric minimum action method in the context of the large deviations theory. By examining the most probable transition paths from metastable to oscillatory regimes for various external carbon input rates, we find two different transition patterns, which gives us an early warning sign for the dramatic change in the carbonate state of the ocean.


\end{abstract}
{\bf Lead paragraph:}

\textbf{
Human activities have been producing more and more carbon dioxide into the carbon cycle, which in turn, significantly influence the climate nowadays. We investigate an oceanic carbonate system that plays a major role in the global carbon cycle. This system has metastable and oscillatory regimes. Due to the small random fluctuations to the external carbon input rate, this system undergoes a transition between these two regimes. We use a geometric minimum action method to capture this transition phenomenon and compute the most probable transition path. We uncover that as the external \ce{CO_2} input rate  $\nu$ changes, two different transition patterns occur with $\nu  \approx 0.2$ as a critical value: (i) At a lower level external \ce{CO_2} input rate (lower than $20 \%$), the concentration of \ce{CO^{2-}_3}  undergoes a larger excursion to shift to the oscillatory state. (ii) But at a higher level external \ce{CO_2} input rate (larger than $20 \%$), a much smaller excursion of the concentration of \ce{CO^{2-}_3} leads to a transition to the oscillatory state. Moreover, as the external \ce{CO_2} input rate $\nu$ increases, the arrival value of the concentration of \ce{CO^{2-}_3} decreases from over 150 $\mathrm{\mu mol \cdot kg^{-1}}$ to around 50 $\mathrm{\mu mol \cdot kg^{-1}}$. 
}

\medskip

\section{Introduction}

The carbon cycle is the biogeochemical cycle in which photosynthesis converts carbon dioxide(\ce{CO_2}) to organic carbon and respiration converts organic carbon back to \ce{CO_2}. Human activities are adding more and more carbon into the atmosphere and become an important part of the carbon cycle. Recently, due to global warming, extreme weather has became more and more common. As a main part of Earth's carbon cycle, the oceanic carbonate system worths our attention. 

Rothman \cite{Rothman14813} introduced a model to describe the dynamical behaviors of the carbonate system of the upper ocean. He took human activities and volcanic emissions as an external source of carbon dioxide (\ce{CO_2}). Using the model, he demonstrated that as the external  \ce{CO_2} inputs increases, both the dissolved inorganic carbon $w$ and carbonate ions \ce{CO^{2-}_3} near a stable state may be disrupted or excited to be near a higher oscillatory level, and eventually come back to the stable state. These large amplitude oscillations are often associated to climate change and mass extinctions.  However, Rothman did not take random fluctuations in the strength of the external source of \ce{CO_2}. But due to unexpected events and measurement error, it is reasonable to take such fluctuations into consideration. 

We look into the same model as Rothman, but include inevitable random fluctuations to the external carbon input. We are interested in the bistable case of this model. Because in original Rothman's model, if the system starts near the stable state, it will never reach the oscillatory state. But in the presence of fluctuations, no matter how small it may be, in a long run, this system will stop staying at the original stable state and travel to an oscillatory state \cite{ditlevsen1999observation,qiu2000kuroshio}. These fluctuations make the Rothman's deterministic model a stochastic dynamical system.

Stochastic dynamical systems are mathematical models for complex phenomena in physical, chemical and biological sciences \cite{arnold2013random,duan2015introduction,duan2014effective,imkeller2012stochastic}. A stochastic system may possess multistable regimes, such as one stable state and one stable limit cycle~(the oscillatory state), with a sandwiched unstable limit cycle. There are many literature concerning the first passage time of transition from the stable fixed point to the limited cycle \cite{li2020particle, li2019first,  li2019influences,li2016levy, wang2016levy}. The transition behaviors are also related to early warning signs \cite{ma2020precursor, ma2019predicting}. We study the transitions between stable regimes under random fluctuations \cite{Creaser2018SequentialNE, PhysRevLett, Lin2019QuasiPotentialCA}. Large deviations theory is a useful tool that can be used to capture transition behaviors of stochastic dynamical system. Roughly speaking, large deviations theory measures how small the probability of rare events is. It generalized the concept of potential of the gradient system to quasi-potential of the non-gradient system, which is a minimizer of the action of paths connecting two ends. We regard such a minimizer as our most probable transition path that connects the metastable state and oscillatory state. We use the geometric minimum action method, which was proposed by M.~Heymann and E.~Vanden-Eijnden \cite{Eric2008}, to find the minimizer. It used a reparametrization method to transfer infinity time scale to finite time scale and then numerically solved an Euler-Lagrange equation to obtain the minimizer.

We will look into the bistable case of the stochastic carbonate system. We will use the geometric minimum action method to find the most probable transition path for different external \ce{CO_2} input rate. Based on the numerical simulations, we uncover two different transition patterns.

This paper is arranged as follows. In section \ref{CCS}, we will briefly introduce the carbon cycle model. The theoretical background of the geometric minimum action method will be provided in section \ref{LDP} and in section \ref{NE}, we will show our numerical experiments and explanation.
 
\section{An oceanic carbonate system}\label{CCS}

\begin{figure}
    \centering
    \includegraphics[scale=0.22]{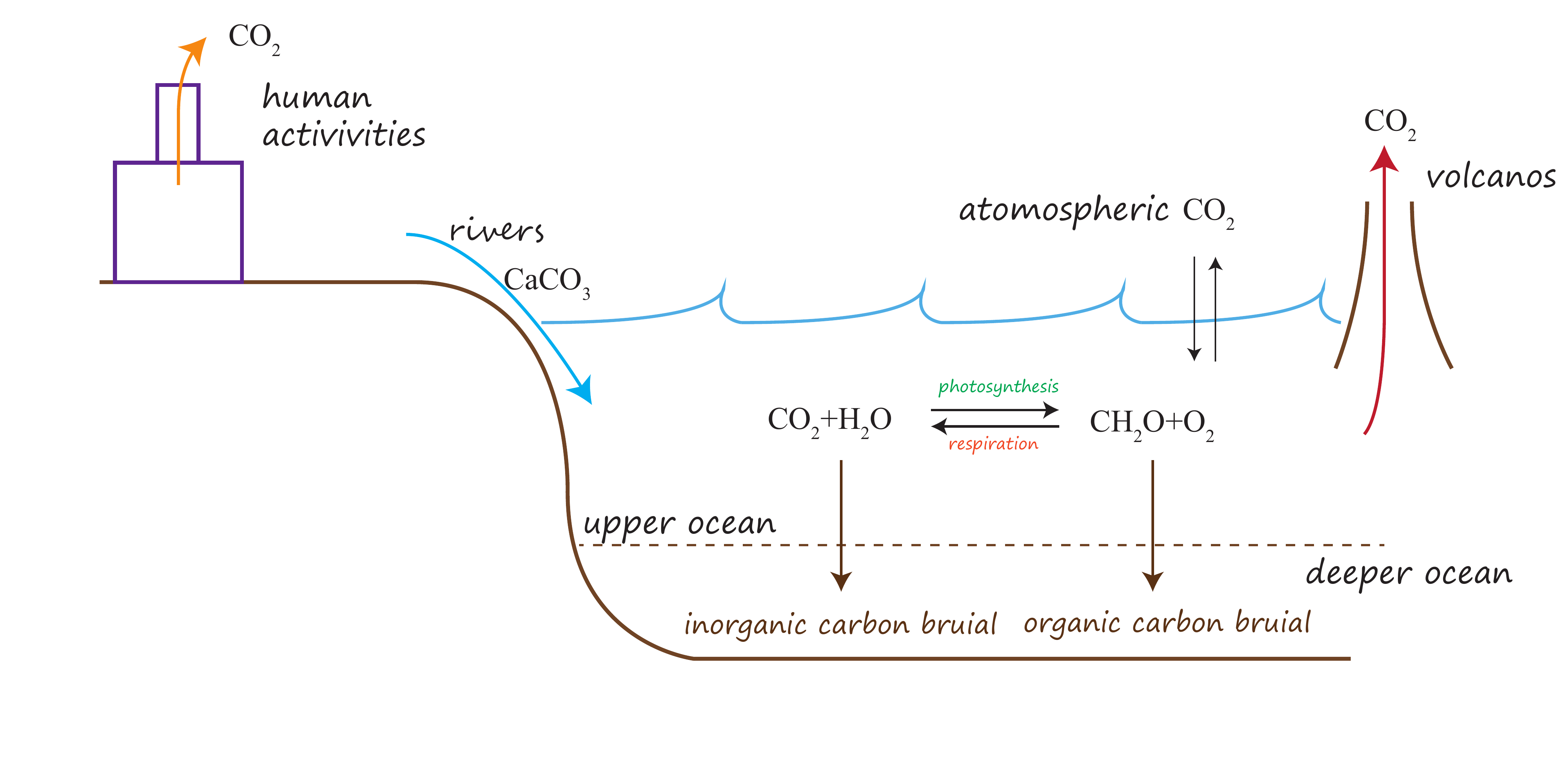}
    \caption{An oceanic carbonate system}\label{OS}
\end{figure}

In the marine carbon cycle, we look into the evolution of the carbonate system in the upper ocean, which is illustrated in Figure \ref{OS}. We consider a well mixed open system and investigate the carbonate ions in the form of \ce{CO^{2-}_3}, \ce{HCO_3^-} and \ce{CO_2}. A main carbonate input is the dissolved \ce{CaCO_3} carried in by rivers. Respiration is a main source of \ce{CO_2}. Volcanic emissions, human activities are regarded as an external source of \ce{CO_2}. The precipitation and transportation to lower ocean of \ce{CaCO_3} is a main carbonate output of the system. Photosynthesis is a main export of \ce{CO_2} of the system. 

Rothman \cite{Rothman14813} introduced a model to describe dynamical behaviors of the carbonate system and it is formulated as follows.
\begin{align}
    dc&=\left[\mu\left[1-b s\left(c, c_p\right)-\theta \bar{s}\left(c, c_x\right)-\nu \right]+w-w_0\right]f(c)dt,\label{CCM-1} \\
    dw&=\left[\mu\left[1-b s\left(c, c_p\right)+\theta \bar{s}\left(c, c_x\right)+\nu\right]-w+w_0\right]dt.\label{CCM-2}
\end{align}
Here, $w$ denotes the total dissolved inorganic carbon, that is 
\begin{equation}
    w=\left[\mathrm{CO}_2\right]+\left[\mathrm{HCO}_3^-\right]+\left[\mathrm{CO}_3^{2-}\right],
\end{equation}
and $c$ is the concentration of \ce{CO^{2-}_3}, which are expressed in units of \ce{\mu mol . kg^{-1}}. The strength of external \ce{CO_2} injection is $\nu$ and is of our main interest in this study. Time $t$ is nondimensionalized by dividing the homeostat’s dominant characteristic timescale $\tau_w$ which is about $10^5$ years. Function $\bar{s}=1-s$. $s$ is a sigmoidal function, and $f$ is the “buffer function” as follow
\begin{equation}
s\left(c, c_{p}\right)=\frac{c^{\gamma}}{c^{\gamma}+c_{p}^{\gamma}},\ \  f\left(c\right)=f_0\frac{c^{\beta}}{c^{\beta}+c_{f}^{\beta}}.
\end{equation}

The parameter $\mu$ is a characteristic concentration, $b$ is the maximum \ce{CaCO_3} burial rate, $\theta$ is the maximum respiration feedback rate, and $\nu$ is the \ce{CO_2} injection rate. The parameter $c_x$ is the crossover \ce{CO^{2-}_{3}}(respiration), $c_p$ is the crossover \ce{CO^{2-}_{3}} (burial), $c_f$ is the crossover \ce{CO^{2-}_{3}} (buffering),   $f_0$ is the maximum buffer factor, and $\gamma,\beta$ are the sigmoid sharpness indexes.
The parameters $\mu$, $b$, $\theta$, $c_p$, $c_x$, $c_f$, $w_0$ $\gamma$ and $\beta$ are constants and are set to fit in properties of the modern ocean.  The value of these parameters are listed in \cite[SI Appendix, Table S1]{Rothman14813}.  We only change the value of $c_x$ and $\nu$ in this article.

\begin{figure}[htbp]
    \centering
    \includegraphics[scale=0.3]{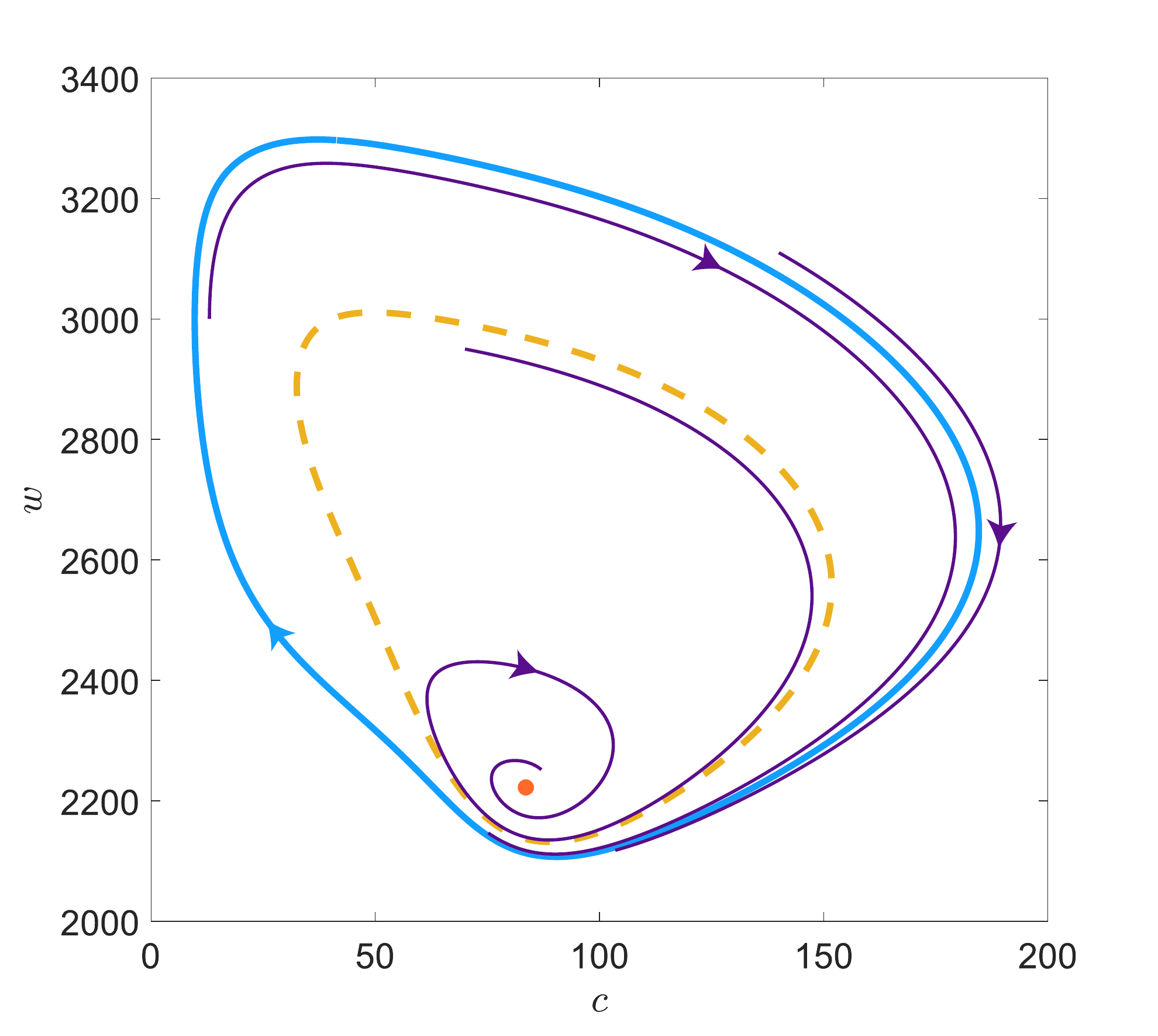}
    \caption{Phase-space trajectories in the bistable regime. Here $c_x=57 \mathrm{ \mu mol \cdot kg^{-1}}$. The yellow dashed limit cycle is unstable. Trajectories (purple arrow lines) initialized inside  the unstable limit cycle return to the stable fixed point (the reddish orange point). Trajectories (purple arrow lines) initialized outside the unstable limit cycle evolve to the stable limit cycle (the blue arrow loop).}
    \label{doublecycle}
\end{figure}

System (\ref{CCM-1}-\ref{CCM-2}) behaves as a Hopf bifurcation when parameter $c_x$ varies. For $c_x<55.89$ $\mathrm{ \mu mol \cdot kg^{-1}}$, there is only one stable fixed point $(c^*,w^*)$ of the system. For $ 55,89$ $ \mathrm{ \mu mol \cdot kg^{-1}}<c_x<62.61$ $ \mathrm{ \mu mol \cdot kg^{-1}} $, a stable fixed point and a stable limit cycle appear and an unstable limit cycle lies in between as is shown in Figure \ref{doublecycle}. For $c_x > 62.61$ $ \mathrm{ \mu mol \cdot kg^{-1}}$, there is a stable limit cycle and an unstable fixed point. For detailed analysis and explanation of the model, we refer readers to \cite{Rothman14813}.

In the real world, there are unexpected events that would slightly change the carbon input rate $\nu$. For example, the Anak Krakatau volcano erupted suddenly in April, 2020 and released massive carbon dioxide  compared to a common volcanic eruption. But this amount of carbon dioxide is still very small compared to the total released carbon dioxide characterized by $\nu$. The instrument error may also lead to disruption to the carbon input rate $\nu$. So, we include a random perturbation to the parameter $\nu$ with sufficient small amplitude $\epsilon$. This gives us the following stochastic system:

\begin{align}
    dc&=\left[\mu\left[1-b s\left(c, c_{p}\right)-\theta \bar{s}\left(c, c_{x}\right)-\nu \right]+w-w_{0}\right]f(c)dt-\epsilon \mu f(c)dB_t^1 \label{CbC-1},\\
    dw&=\left[\mu\left[1-b s\left(c, c_{p}\right)+\theta \bar{s}\left(c, c_{x}\right)+\nu\right]-w+w_{0}\right]dt+\epsilon \mu dB_t^2, \label{CbC-2}
\end{align}
where $B_t^1$ and $B_t^2$ are two standard Brownian motions.

No matter how small the random perturbation is, the original stable state would become a metastable state. And after a sufficient long time, the solution to (\ref{CbC-1}-\ref{CbC-2}) would transport from one metastable state to another metastable state. That means, for the parameter $c_x$ lying in the bistable regime $ (55,89$ $ \mathrm{ \mu mol \cdot kg^{-1}},62.61$ $ \mathrm{ \mu mol \cdot kg^{-1}} )$, the solution to system (\ref{CbC-1}-\ref{CbC-2}) will go from the stable state $(c^*,w^*)$ to the oscillatory state, that is the stable limit cycle of the deterministic system (\ref{CCM-1}-\ref{CCM-2}).

We will focus on the bistable case, that is $ 55,89$ $ \mathrm{ \mu mol \cdot kg^{-1}}<c_x<62.61$ $ \mathrm{ \mu mol \cdot kg^{-1}} $. In the next section, we will introduce a large deviation approach that enables us to capture this transition phenomena.

\section{A large deviation approach}\label{LDP}
We will briefly introduce the large deviations theory and the geometric minimum action method based on it.

We consider a generalized stochastic differential equation in $\mathbb{R}^2$ of equation (\ref{CbC-1}-\ref{CbC-2}) as follows,
\begin{equation}\label{SDE}
dX^\epsilon (t)=\kappa  (X^\epsilon(t))dt+\epsilon \eta  (X^\epsilon(t))dB(t), \quad X^\epsilon(0)=x_0,
\end{equation}
where $\kappa:\mathbb{R}^2 \rightarrow \mathbb{R}^2$ is a regular function, $\eta:\mathbb{R}^2\rightarrow\mathbb{R}^{2\times 2}$ is a $2\times 2$ matrix-valued function and $B$ is a standard Brownian motion in $\mathbb{R}^2$. 

Under some assumptions on the coefficients $\kappa$ and $\eta $, for instance, the Lipschitz continuous conditions, the solution $X^\epsilon $ to equation (\ref{SDE}) satisfies a large deviation principle with an action functional $S_T$. It enables us to estimate the probability that $X^\epsilon$ stays in a $\delta$-neighborhood of a path $\psi$ though the following asymptotic relation.

\begin{equation}
    \mathbb{P}(d(X^\epsilon,\psi)<\delta) \sim \exp \left(-\frac{1}{\epsilon} S_T(\psi)\right), \quad \textrm{as } \epsilon \to 0,
\end{equation} 
where $d$ is the distance in $\mathbb{R}^2$.
The action functional $S_T$ is
\begin{equation}
S_T(\varphi)=\int_0^T L(\varphi,\dot{\varphi})dt,
\end{equation}
with the Lagrange functional  
\begin{equation*}
L(\varphi,\dot{\varphi})=\frac{1}{2}\left(\dot{\varphi}-\kappa(\varphi)\right)(\eta\eta^T(\varphi))^{-1}[\dot{\varphi}-\kappa(\varphi)].
\end{equation*}
The large deviation theory also equips us to capture the long time behaviors of $X^\epsilon$ with quasi-potential $V$. That is 
\begin{equation}\label{qpot}
    V(x_1,x_2)=\inf_{T>0}\inf_{\varphi \in \bar{C}_{x_1}^{x_2}(0,T)} S_T(\varphi),
\end{equation}
where $\bar{C}_{x_1}^{x_2}(0,T)$ is the space of all absolutely continuous functions that start at $x_2$ and end at $x_2$. Assume that there exists a path $\tilde{\varphi}$ and time $\tilde{T} \in [0,\infty]$ satisfying,
\begin{equation*}
    V(x_1,x_2)=S_{\tilde{T}}(\tilde{\varphi}).
\end{equation*}
Then, in probability $1$, the sample paths of $X^\epsilon$ from $x_1$ to a $\delta$-neighborhood of $x_2$, converge to $\tilde{\varphi}$ when the noise intensity $\epsilon$ and the neighborhood size $\delta$ tend to $0$ as well as the time $T$ tends to $\tilde{T}$. This result means that the sample path $\tilde{\varphi}$ is the most probable transition path we are looking for. Readers may refer to \cite[Proposition 2.3]{Eric2008} for a more precise statement.

In summary, we only need to find the minimizer of equation (\ref{qpot}) to obtain a most probable transition path between two metastable states. For a most probable transition path between a metastable state and an oscillatory state $D$, for example, a limit cycle, we need to find the minimizer of the following equation.
\begin{equation} \label{mini}
    V(x,D)=\inf_{y \in D} V(x,y)=\inf_{y \in D} \inf_{T>0}\inf_{\varphi \in \bar{C}_{x}^{y}(0,T)} S_T(\varphi).
\end{equation}

We will conduct a geometric minimum action method introduced in \cite{Eric2008} to find the most probable transition path between the stable fixed point and points near the limit cycle. Then, we pick up the path which has the minimum quasi-potential as the most probable transition path between the stable state and the limit cycle. That is the minimizer to equation (\ref{mini}).

\section{Numerical simulations and discussions}\label{NE}
In this section, we will present the most probable transition paths that are computed through the geometric minimum action method. We are interested in how the most probable transition path changes as the external \ce{CO_2} input rate $\nu$ varies, in order to better understand this carbonate system. 

\begin{figure}[htbp]
    \centering

    \subfloat[$\nu = 0$ ]{
        \label{nu0}
        \includegraphics[width=0.48\textwidth]{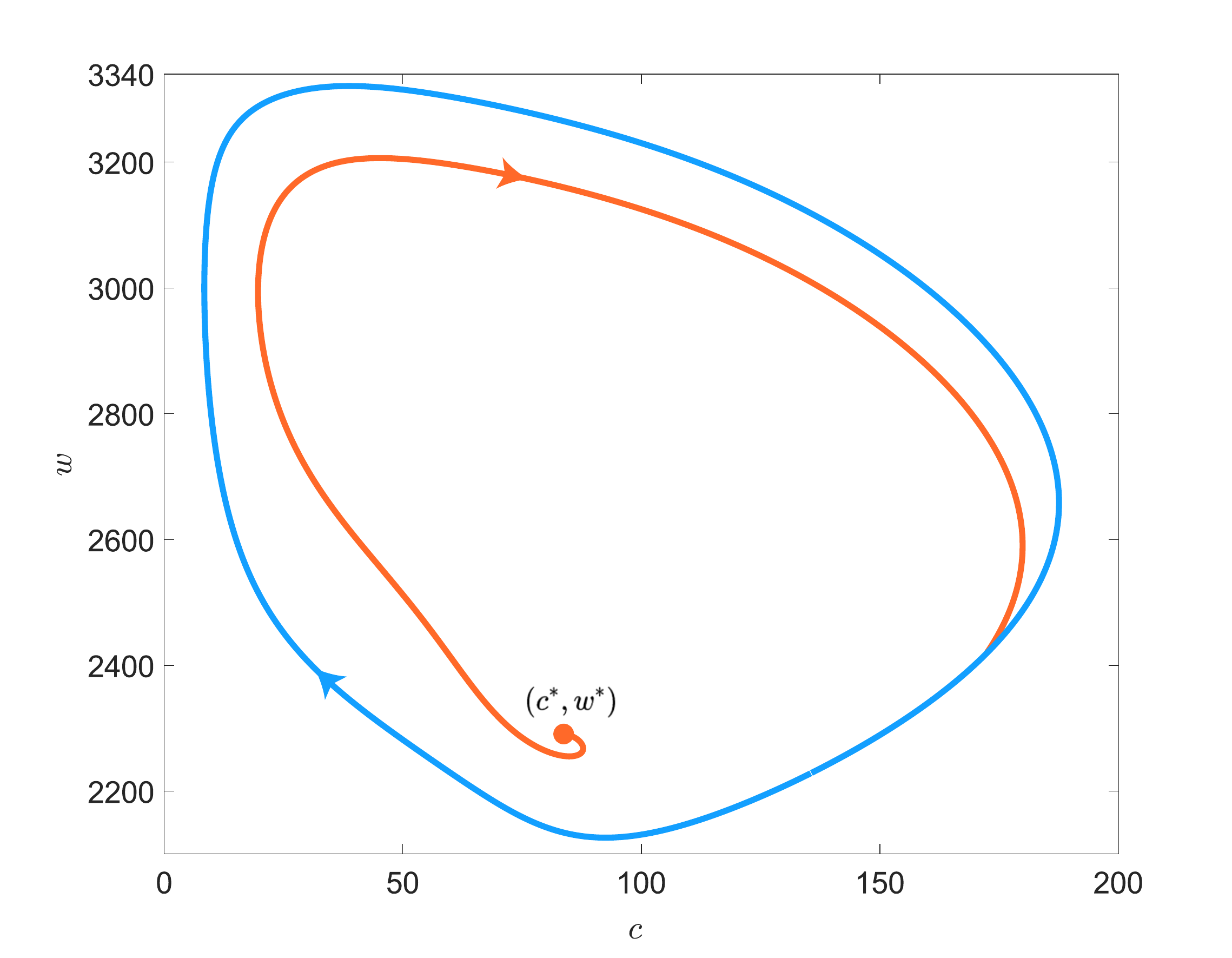}
    }
    \subfloat[$\nu = 0.19$ ]{
        \label{nu019}
        \includegraphics[width=0.48\textwidth]{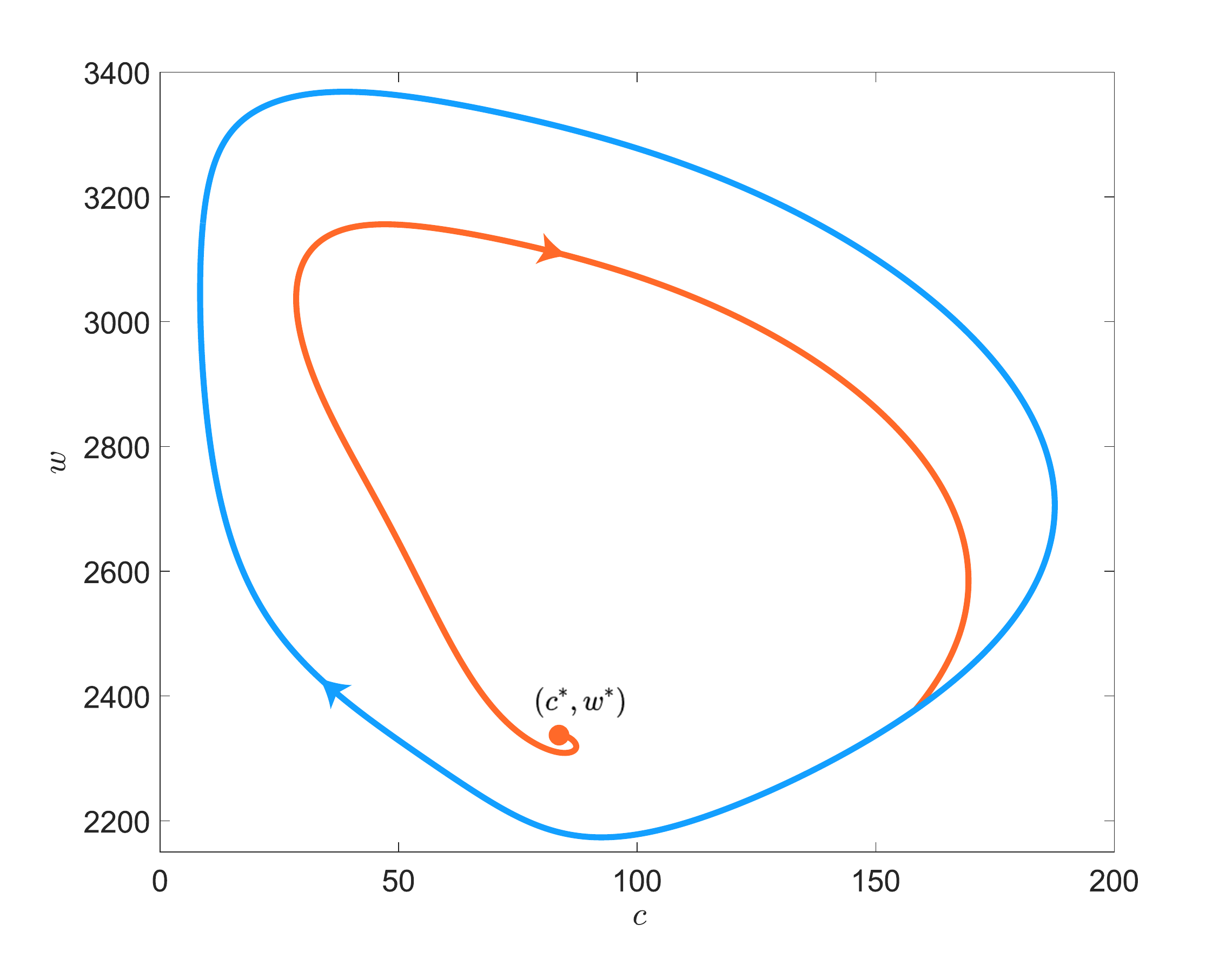}
    }
\\
    \subfloat[ $\nu = 0.2$]{
        \label{nu04}
        \includegraphics[width=0.48\textwidth]{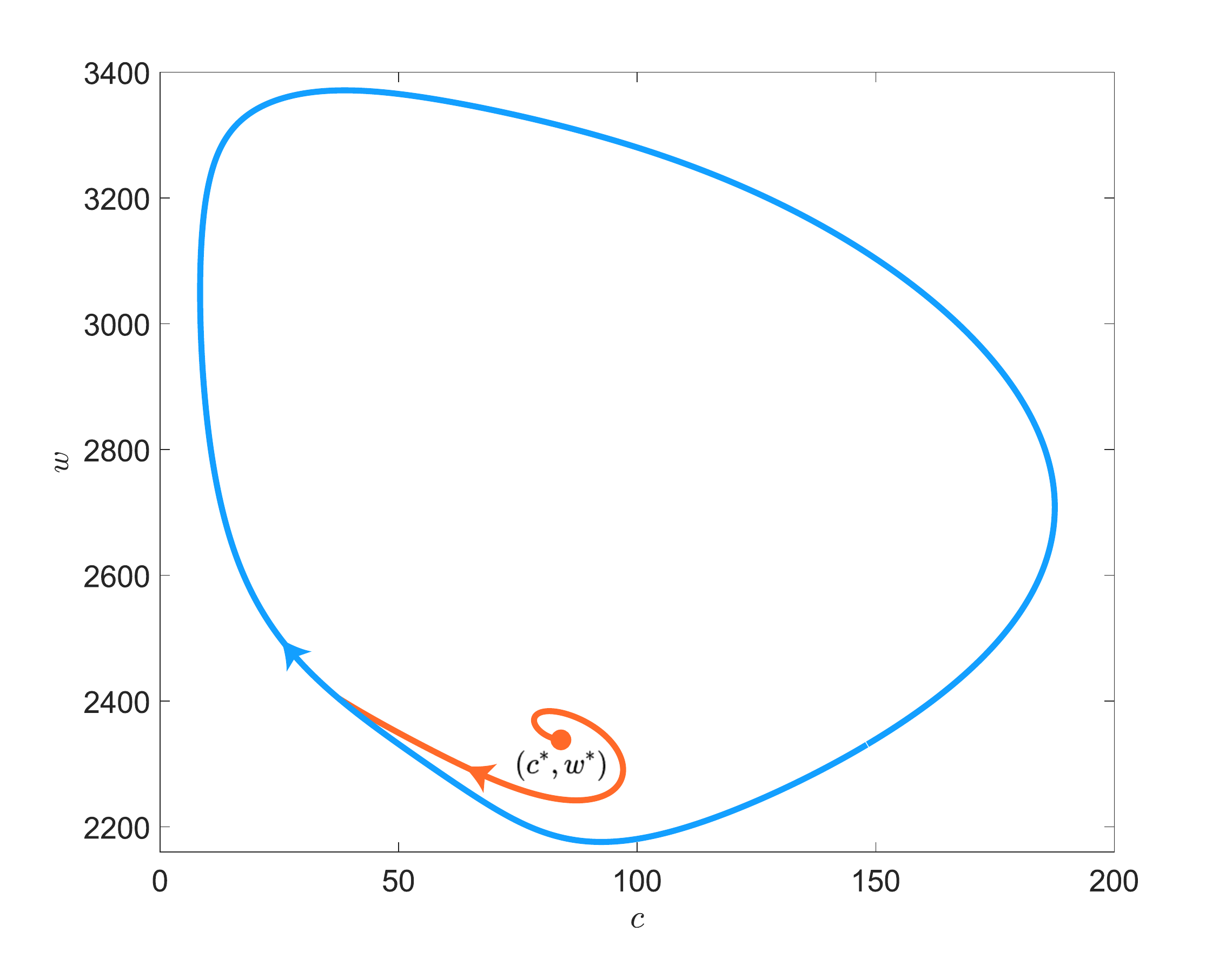}
    }
    \subfloat[$\nu = 0.4$ ]{
        \label{nu09}
        \includegraphics[width=0.48\textwidth]{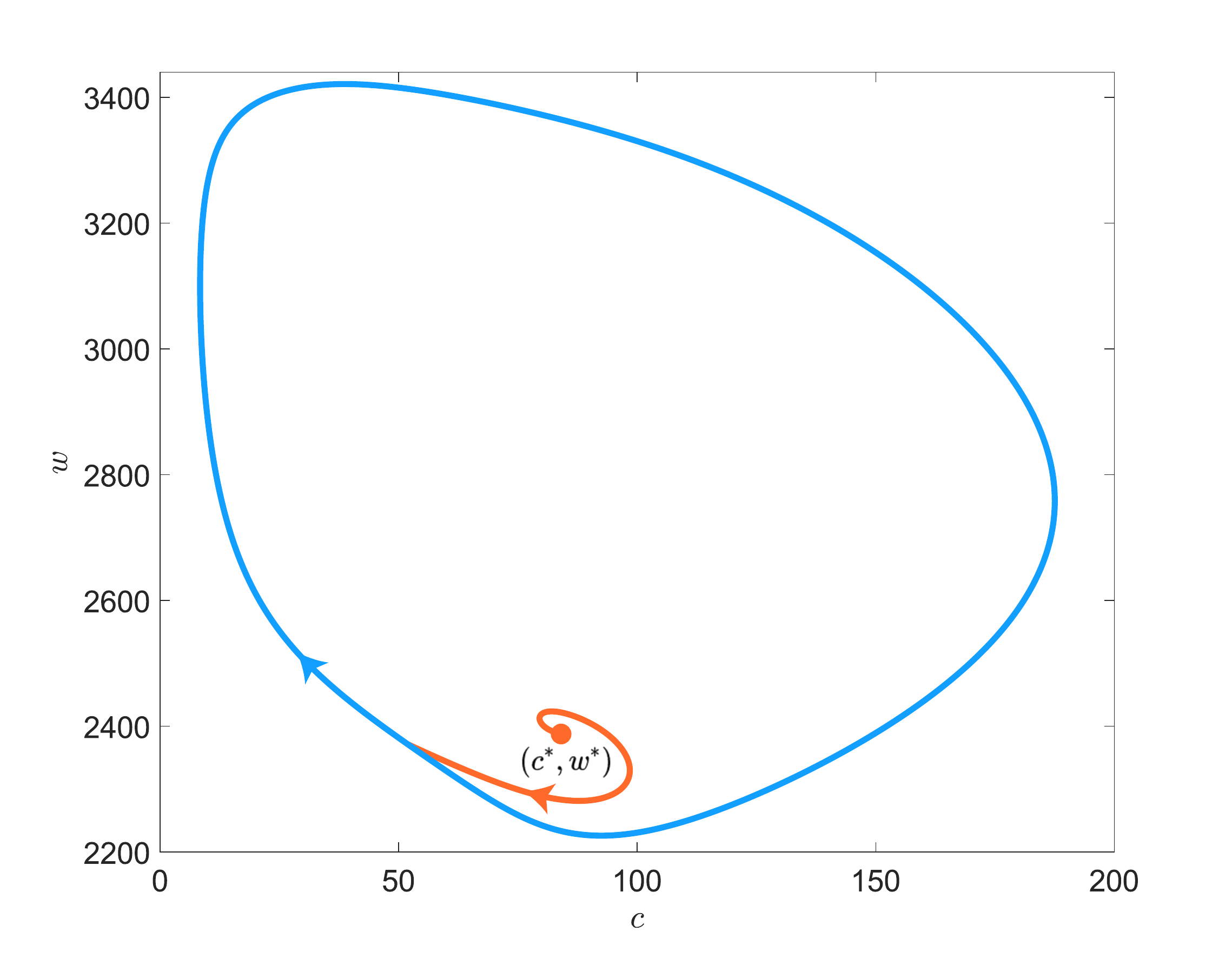}
    }
    \caption{Most probable transition paths~( reddish orange arrow lines) from the metastable state $(c^*,w^*)$ (orange point) to the oscillatory state (blue limit cycle) for external \ce{CO_2} input rate $\nu=0, 0.19, 0.2 0.4$  and parameter $c_x=62$. }\label{tras1}
\end{figure}

We choose $c_x=62$ and use the Euler scheme to simulate the limit cycle value, that is the oscillatory state of the carbonate system. We calculate the most probable transition path from the stable state $(c^*,w^*)$ to a point near the oscillatory state, for $\nu$ from $0$ to $0.9$. $3000$ points on the transition path are made to be equidistant. The quasi-potential between two points on the limit cycle is zero, because there is a deterministic orbit connecting each other though the limit cycle and this makes the action achieve its minimum $0$. In practice, as we do not actually choose points on the limit cycle, we can obtain a most probable transition path though the geometric minimum action method. Then, we only need to find the path with the least quasi-potential as the most probable transition path connecting the stable state $(c^*,w^*)$ and the oscillatory state. The most probable transition paths for $\nu=0,0.19,0.2$ and $0.4$ are shown in Figure \ref{tras1} in the $c$-$w$ plane.

\begin{figure}[htbp]
    \centering
    \subfloat[The length  for the most probable transition path]{
        \includegraphics[width=0.48\textwidth]{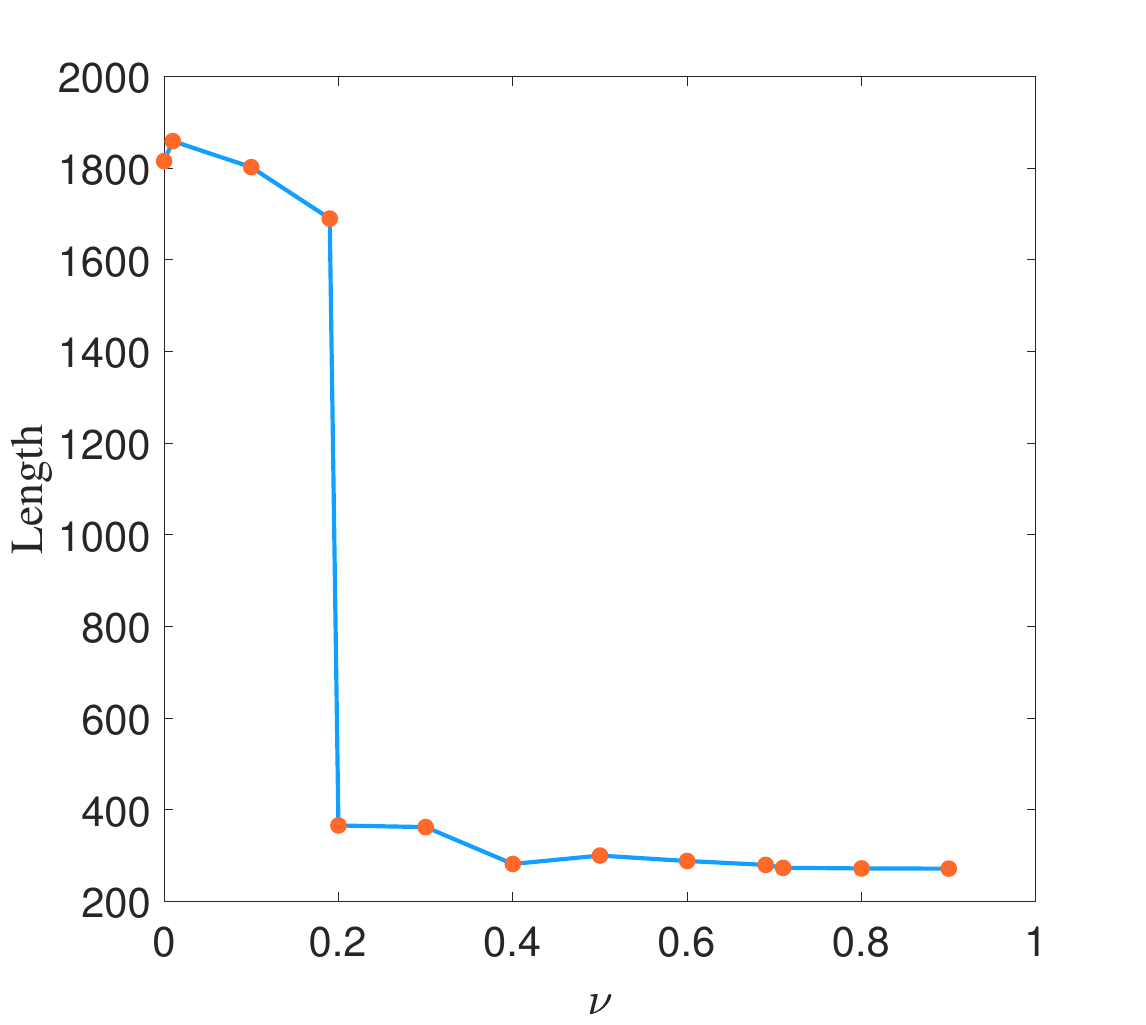}
        \label{nulen}
    }
    \subfloat[The arrival value of $c$ ]{
        \label{nuend}
        \includegraphics[width=0.48\textwidth]{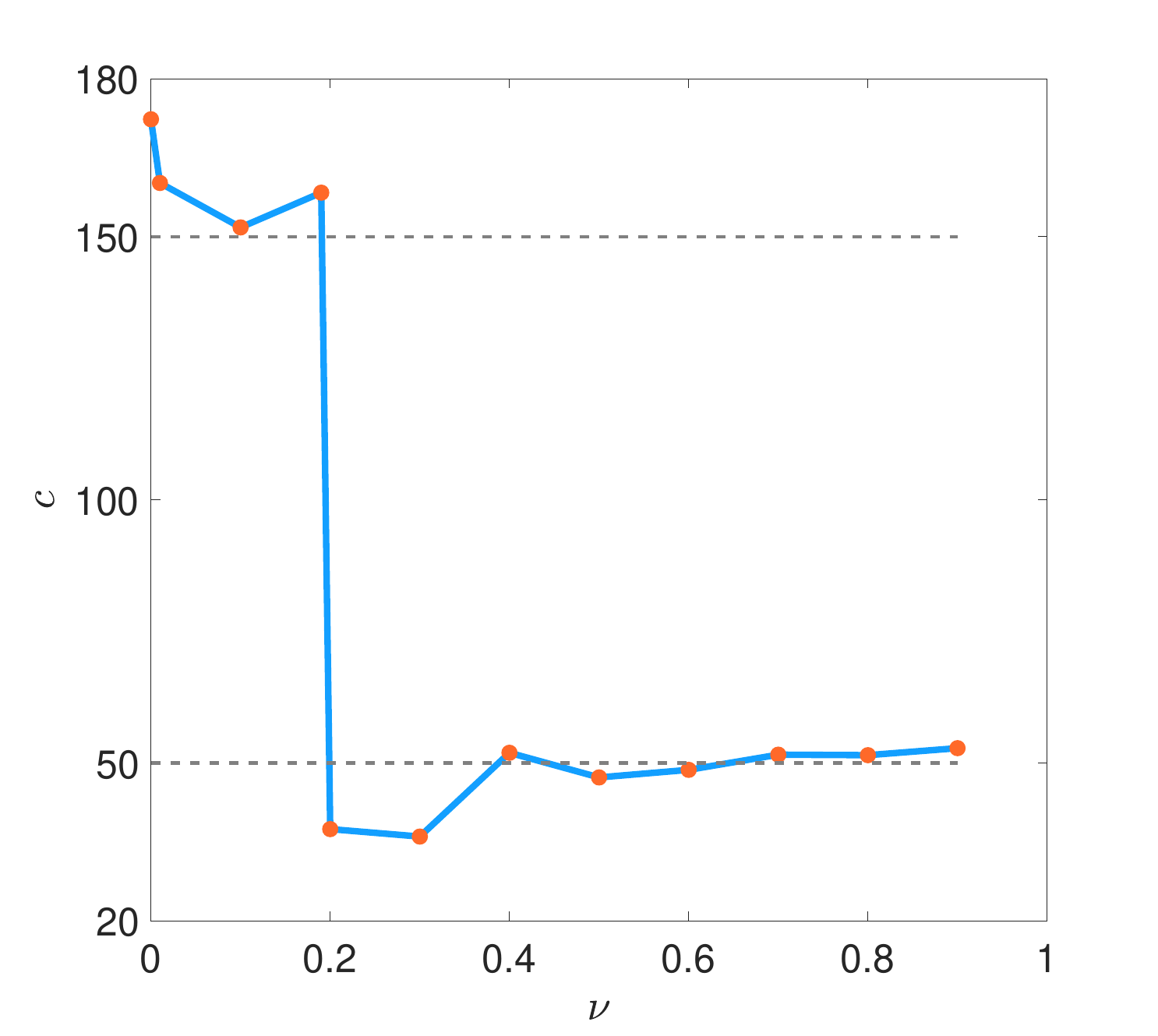}
    }
    \caption{The length for the most probable transition path and the arrival value of concentration of \ce{CO^{2-}_3}, $c$, with respect to different external \ce{CO_2} input rate $\nu$. Parameter $c_x=62.$ Critical changes happen around $\nu=0.2$. The length of most probable transition path decreases to lower than 400. The arrival value of $c$ decreases from 150 $\mathrm{\mu mol \cdot kg^{-1}}$ to around 50 $\mathrm{\mu mol \cdot kg^{-1}}$.}\label{tras3}   
\end{figure}

The most probable transition path has to travel a longer way to achieve the oscillatory state for $\nu < 0.2$ as is shown in Figure \ref{tras3}. The dramatic change happens near $\nu=0.2$. The most probable transition path stops ending at the ``southeast" part of the limit cycle and starts to end at the ``southwest'' part of the limit cycle as is shown in Figure \ref{tras1}. Consequently, the oscillatory ranges for $c$ and $w$ on the most probable transition path shrink. And the arrival value of concentration of \ce{CO^{2-}_3} decreases from above 150 $\mathrm{\mu mol \cdot kg^{-1}}$ to around 50 $\mathrm{\mu mol \cdot kg^{-1}}$.

\begin{figure} 
    \centering
    \subfloat[ Transition in view of $c$ for $\nu=0.19$]{
        \label{cnu01}
        \includegraphics[width=0.48\textwidth]{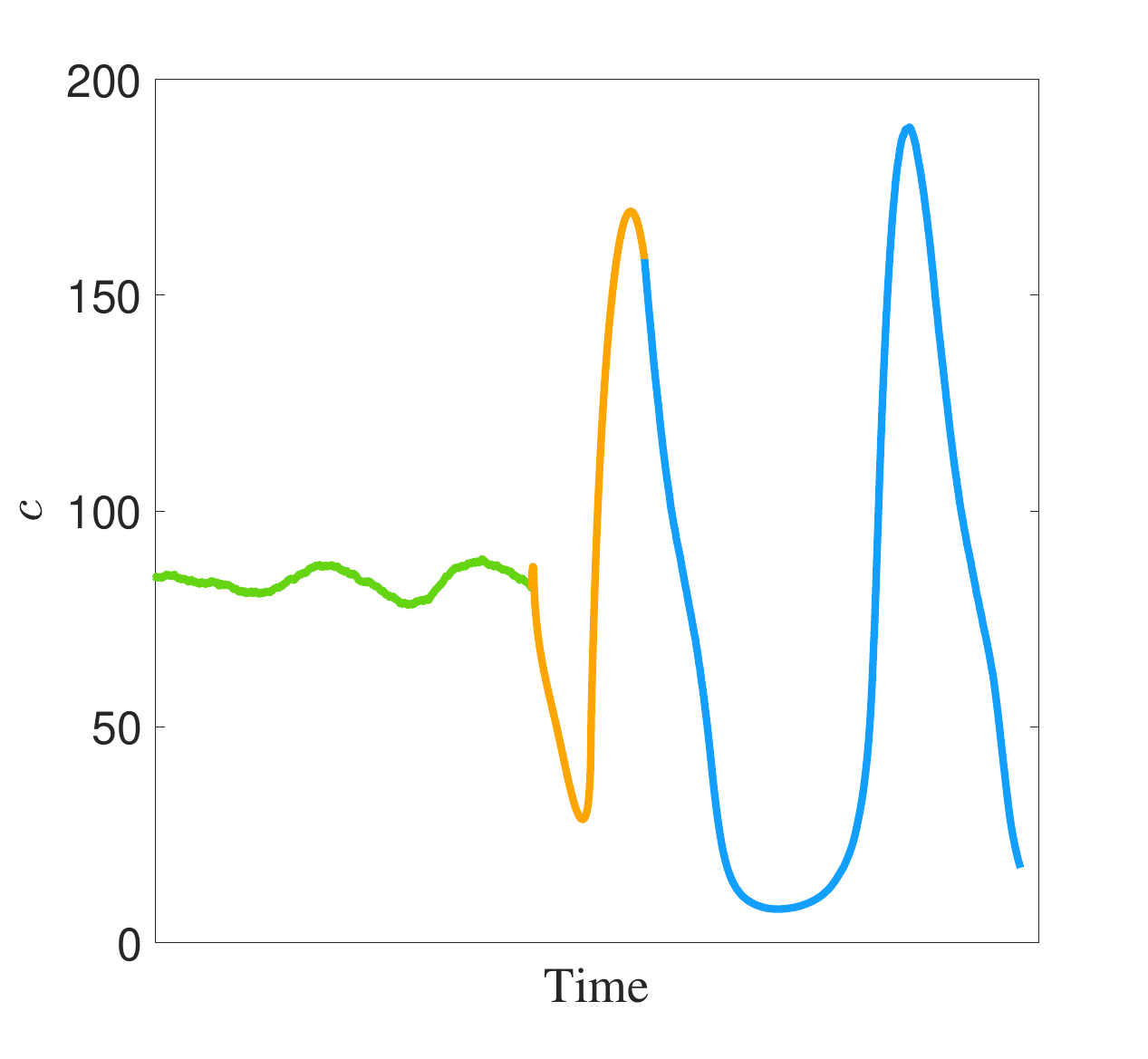}
    }
    \subfloat[ Transition in view of $c$ for $\nu=0.2$]{
        \label{cnu04}
        \includegraphics[width=0.48\textwidth]{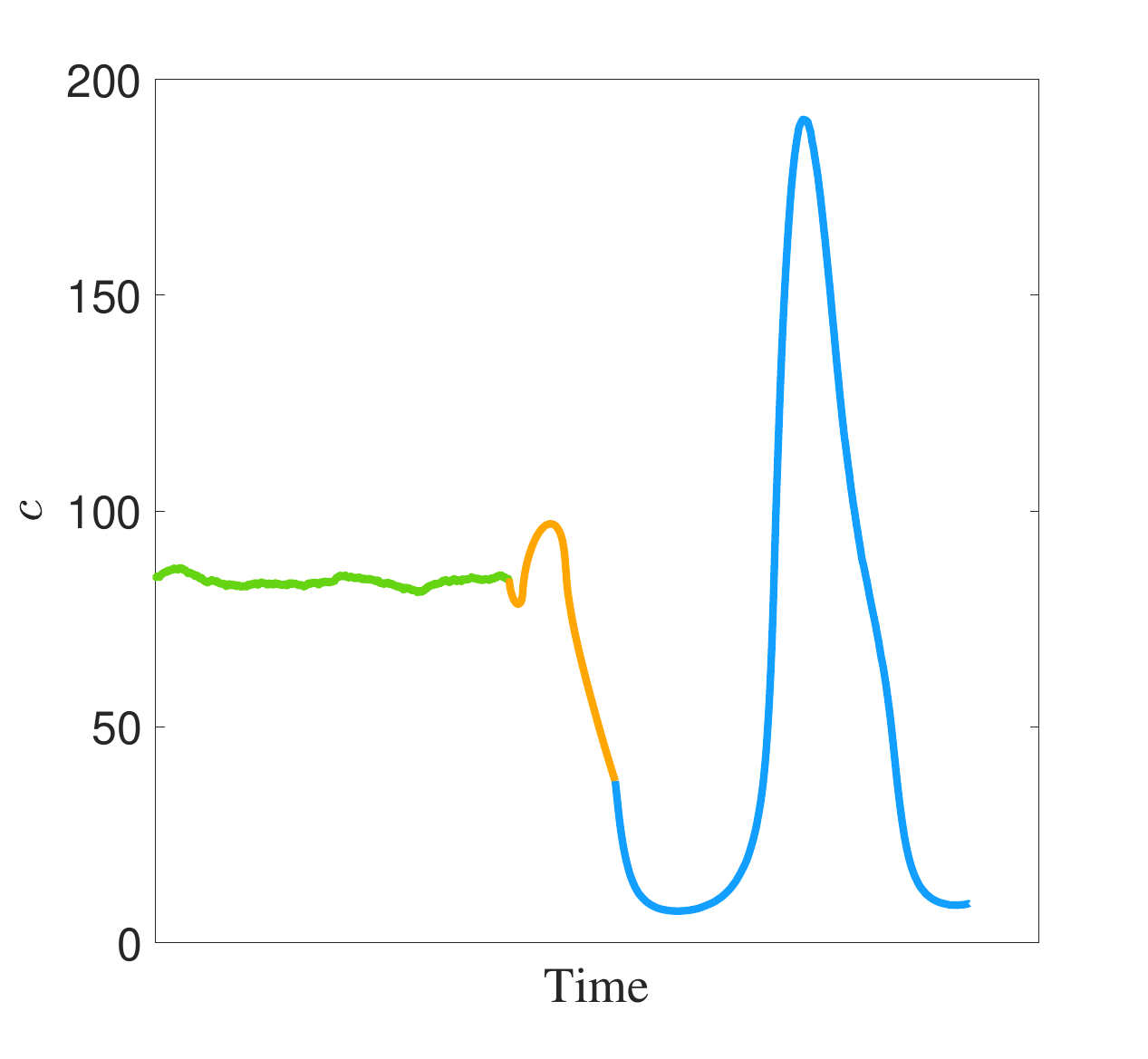}
    }
    \caption{Illustrations for transitions in view of $c$, the concentration of \ce{CO^{2-}_3}, for \ce{CO_2} input rate $\nu=0.19, 0.2$  and parameter $c_x=62$. The orange line is the value of $c$ on the most probable transition path. The green line represents the value of $c$ near the metastable state and the blue line represents the value of $c$ on the oscillatory state.}\label{twopatterns}
\end{figure}
If we look into how $c$, the concentration of \ce{CO^{2-}_3}, changes with respect to time, we find two different patterns for transitions. Again, we use Euler scheme to generate the data of $c$ before and after transition for noise strength $\epsilon = 0.01$. Then we connect them with the most probable transition path and illustrate the results in in Figure \ref{tras3}. When the the external \ce{CO_2} rate $\nu$ is smaller than a threshold near $0.2$, the concentration of \ce{CO^{2-}_3} will undergo an excursion, whose size is slightly smaller than the oscillatory state (the limit cycle shown in Figure \ref{tras1}, to shift to the oscillatory state. When $\nu$ exceeds the threshold near $0.2$, the concentration of \ce{CO^{2-}_3} will undergo a much smaller excursion to shift to the oscillatory state. This result implies that at more that $20 \%$ external carbon input rate, much less changes to the concentration of \ce{CO^{2-}_3} may lead to a transition. 

Our calculation shows the early warning sign of the carbonate system. If the solution $(c,w)$ stays near the stable state $(c^*,w^*)$, the carbonate system of the ocean and the carbon cycle are in dynamic equilibrium. But when the solution $(c,w)$ enters the oscillatory state (the stable limit cycle), the carbonate system of the ocean will change dramatically and may cause massive extinction \cite{Rothman14813}. It is important to identify the transition between these two states. Our calculation shows two transition patterns of $c$, the concentration of \ce{CO^{2-}_3}(See Figure \ref{twopatterns}): (i) At a lower level external \ce{CO_2} input rate (lower than $20 \%$), $c$ will first drop to lower than 50 $\mathrm{\mu mol \cdot kg^{-1}}$ and then increase to higher than 150 $\mathrm{\mu mol \cdot kg^{-1}}$, before transferring to the oscillatory state. (ii) But at a higher level external \ce{CO_2} input rate (larger than $20 \%$), $c$ will only increase a little to around to 100 $\mathrm{\mu mol \cdot kg^{-1}}$ and then drop to around 50 $\mathrm{\mu mol \cdot kg^{-1}}$ to enter the oscillatory state. Given the external \ce{CO_2} input rate, we can identify the transition before the system enters the oscillatory state.

%need not to change too much to transfer to another stable state at high level external \ce{CO_2} input rate.

In conclusion, we have computed the most probable transition path connecting a metastable state and an oscillatory state for a carbonate system under random fluctuations in the external \ce{CO_2} input rate. As   the external \ce{CO_2} input rate $\nu$ changes, two different transition patterns occur with $\nu  \approx 0.2$ as a critical value (See Figure \ref{twopatterns}): (i) At a lower level external \ce{CO_2} input rate (lower than $20 \%$), the concentration of \ce{CO^{2-}_3}  undergoes a larger excursion to shift to the oscillatory state. (ii) But at a higher level external \ce{CO_2} input rate (larger than $20 \%$), a much smaller excursion of the concentration of \ce{CO^{2-}_3} leads to a transition to the oscillatory state. Moreover, as the external \ce{CO_2} input rate $\nu$ increases, the arrival value of the concentration of \ce{CO^{2-}_3} decreases from over 150 $\mathrm{\mu mol \cdot kg^{-1}}$ to around 50 $\mathrm{\mu mol \cdot kg^{-1}}$.(Figure \ref{nuend}).

\section*{Acknowledgements}
This work was partly supported by the NSFC grants 11771449.

\section*{Data Availability}
The data that support the findings of this study are openly available in GitHub. \\ \url{https://github.com/JayWeiess/CarbonCycle-gMAM}

\bibliographystyle{abbrv}
 %\bibliographystyle{plain}
% \biboptions{square,numbers,sort&compress}
\bibliography{Manuscript}
 
\end{document}